\documentclass{article}
\usepackage{amssymb,amsmath,multirow}
\usepackage[dvips]{graphicx}
\newtheorem{thm}{Theorem}
\newtheorem{prop}{Proposition}
\newtheorem{lem}{Lemma}
\newtheorem{cor}{Corollary}
\newtheorem{exam}{Example}
\def\pn{\par\smallskip\noindent}
\def\proof{\pn {Proof.} }
\def\endproof{\hfill \quad{$\Box$}\smallskip}

\title{Trust Region Subproblem with a Fixed Number of Additional Linear Inequality Constraints has Polynomial Complexity\thanks{This research was supported by Taiwan National Science Council under grant 102-2115-M-
006-010, by National Center for Theoretical Sciences (South), by National Natural
Science Foundation of China under grant 11001006 and
91130019/A011702, and by the fund of State Key Laboratory of
Software Development Environment under grant SKLSDE-2013ZX-13.}}

\author{Yong  Hsia\footnotemark[2] \and Ruey-Lin Sheu\footnotemark[3]~\footnotemark[4]}

\begin{document}
\maketitle
\renewcommand{\thefootnote}{\fnsymbol{footnote}}

\footnotetext[2]{State Key Laboratory of Software Development
              Environment, LMIB of the Ministry of Education, School of
Mathematics and System Sciences, Beihang University, Beijing,
100191, P. R. China ({\tt dearyxia@gmail.com}).}

\footnotetext[3]{Department of Mathematics, National Cheng Kung
University, Taiwan ({\tt
rsheu@mail.ncku.edu.tw}).}

\footnotetext[4]{Corresponding author.}

\begin{abstract}
The trust region subproblem with a fixed number $m$ additional
linear inequality constraints, denoted by (${\rm T_m}$), have drawn
much attention recently. The question as to whether Problem (${\rm
T_m}$) is in Class P or Class NP remains open. So far, the only
affirmative general result is that (${\rm T_1}$) has an exact
SOCP/SDP reformulation and thus is polynomially solvable. By
adopting an early result of Mart\'{\i}nez on local non-global
minimum of the trust region subproblem, we can inductively reduce
any instance in (${\rm T_m}$) to a sequence of trust region
subproblems (${\rm T_0}$). Although the total number of (${\rm
T_0}$) to be solved takes an exponential order of $m$, the reduction
scheme still provides an argument that the class (${\rm T_m}$) has
polynomial complexity for each fixed $m$. In contrast, we show by a
simple example that, solving the class of extended trust region
subproblems which contains more linear inequality constraints than
the problem dimension; or the class of instances consisting of an
arbitrarily number of linear constraints, namely
$\bigcup_{m=1}^\infty ({\rm T_{m}}$), is NP-hard. When $m$ is small
such as $m=1,2$, our inductive algorithm should be more efficient
than the SOCP/SDP reformulation since at most 2 or 5 subproblems of
(${\rm T_0}$), respectively, are to be handled. In the end of the
paper, we improve a very recent dimension condition by Jeyakumar and
Li under which (${\rm T_m}$) admits an exact SDP relaxation.
Examples show that such an improvement can be strict indeed.
\end{abstract}


\pagestyle{myheadings}
\thispagestyle{plain}

\pagestyle{myheadings} \thispagestyle{plain} \markboth{Y. Hsia, R.L.
Sheu} {Complexity of Extended Trust Region Subproblem}

\section{Introduction}
\label{intro}

The {\it classical} trust region subproblem, which minimizes a
nonconvex quadratic function over the unit ball
\begin{eqnarray}
(\rm T_0)~~~&\min & f(x)=\frac{1}{2}x^TQx+c^Tx \nonumber\\
&{\rm s.t.}&x^Tx\le 1,\nonumber
\end{eqnarray}
is an important feature in trust region methods \cite{C00,Y90}. It is
well known that finding an $\epsilon$-optimal solution of ($\rm
T_0$) has polynomial complexity \cite{F98,Y92} and efficient
algorithms for solving ($\rm T_0$) can be found in
\cite{G99,M83,R97}. Moreover, problem ($\rm T_0$) is also a special
case of a quadratic problem subject to a quadratic inequality
constraint (QP1QC). It was proved that, under Slater's condition,
(QP1QC) admits a tight SDP relaxation and its optimal solution can
be found through a matrix rank one decomposition procedure
\cite{S03,D}.

Extensions of ($\rm T_0$), sometimes termed as the {\it extended}
trust region subproblem, consider problems such as adding to {\rm
($\rm T_0$)} several linear inequality constraints \cite{Y03} or imposing a
full-dimensional ellipsoid \cite{C85}. In particular, we are
interested in the following variant:
\begin{eqnarray}
(\rm T_m)~~~&\min & f(x)=\frac{1}{2}x^TQx+c^Tx \nonumber\\
&{\rm s.t.}&x^Tx\le 1,\label{Tm2}\\
&&a_i^Tx\le b_i,~i=1,\ldots,m,\label{Tm3}
\end{eqnarray}
which arises from applying  trust region methods to solve constrained nonlinear programs  \cite{C00}. We notice that some NP-hard combinatorial
optimization problems also have the similar formulation.
A typical example is to rewrite
the standard quadratic program
\begin{eqnarray*}
(\rm{QPS})~~~&\min& x^TQx\\
&{\rm s.t.}&e^Tx=1,\\
&&x\ge 0
\end{eqnarray*}
as a special case of (${\rm T_m}$), where $e$ is the vector of all ones. To do so, let
$y=(x_1,\ldots,x_{n-1})^T$. By replacing $e^Tx=1,~x\ge0$ with
$x_n=1-e^Ty\ge 0,~y\ge 0$, we can express the standard quadratic
program (QPS) in terms of variable $y$ as
\begin{eqnarray}
&\min& \left(\begin{array}{c}y\\1-e^Ty\end{array}\right)^TQ
\left(\begin{array}{c}y\\1-e^Ty\end{array}\right)\nonumber\\
&{\rm s.t.}&1-e^Ty\ge 0,\label{QPS-reformulation1}\\
&&y\ge 0.\label{QPS-reformulation1a}
\end{eqnarray}
It is easy to see that $0\le y\le e$ and
\[
y^Ty=\sum_{i=1}^{n-1}y_i^2\le \sum_{i=1}^{n-1}y_i\le 1.
\]
In other words, by imposing a redundant constraint $y^Ty\le1$ to
(\ref{QPS-reformulation1})-(\ref{QPS-reformulation1a}), we enforce
(QPS) to have an equivalent extended trust region subproblem
reformulation as follows:
\begin{eqnarray*}
(\rm{QPS-TRS})~~~&\min& \left(\begin{array}{c}y\\1-e^Ty\end{array}\right)^TQ
\left(\begin{array}{c}y\\1-e^Ty\end{array}\right)\\
&{\rm s.t.}&1-e^Ty\ge 0,\nonumber\\
&&y\ge 0,\nonumber\\
&&y^Ty\le 1.\nonumber
\end{eqnarray*}
Since (QPS) is NP-hard (as it captures the NP-hard combinatorial
problem to find the cardinality number of the maximum stable set in
a graph), so is ({QPS-TRS}). Let (${\rm T_{n+1}}$) represent the
class of extended trust region subproblems which always has the
number of linear inequality constraints exceeding the problem
dimension by one. We can immediately conclude from the example
({QPS-TRS}) that (${\rm T_{n+1}}$) must be NP-hard. The implication
is that solving the subclass of extended trust region subproblems
which contains more linear inequality constraints than the problem
dimension; or solving the most general extended trust region
subproblems consisting of an arbitrarily number of linear
constraints, namely $\bigcup_{m=1}^\infty ({\rm T_{m}})$, should be
difficult.

A natural question arises from computational complexity: ``Fix a
positive integer $m$. What is the complexity of solving $({\rm
T_{m}})$ for all possible dimensions?'' The problem turns out to be
more difficult than most people thought. The only affirmative result
so far in literature is that (${\rm T_{m}}$) with ${\rm {m}= 1}$ is
polynomially solvable \cite{B13,S03}. For ${\rm {m}= 2}$, the
polynomial solvability of some special cases of ${\rm (T_2)}$ were
established when $a_1$ and $a_2$ are parallel \cite{B13,Y03}; or
when $a_1^Tx\le b_1$ and $a_2^Tx\le b_2$ are non-intersecting in the
unit-ball \cite{B14}. When $m\ge 2$ and any two inequalities are
non-intersecting in the interior of the unit ball, this subclass of
(${\rm T_m}$) is also polynomial solvable as shown in \cite{B14}.
Very recently, Jeyakumar and Li \cite{J13} showed that, under the
following {\bf dimension condition}, (${\rm T_{m}}$) is also
polynomial solvable \cite{J13}:
\begin{equation}
[{\rm DC}]~~{\rm dim~Ker}(Q-\lambda_{\min}(Q)I_n)\ge {\rm dim~span}\{a_1,\ldots,a_m\}+1,\label{dc}
\end{equation}
where Ker$(Q)$ denotes the kernel of $Q$; $\lambda_{\min}(Q)$ the
smallest eigenvalue of $Q$; dim $L$ the dimension of a subspace $L$;
and $I_n$ the identity matrix of order $n$.

All the approaches mentioned above elaborate the polynomial
complexity of some (${\rm T_{m}}$) through an exact SOCP/SDP
reformulation \cite{B13,B14,S03} such as, for $m=1$,
\begin{eqnarray}
 &\min& \frac{1}{2}{\rm trace}(QX)+c^Tx \nonumber \\
&{\rm s.t.}&  {\rm trace}(X)\le 1,~X\succeq xx^T,\nonumber\\
&&\|b_1x-Xa_1\|_2\le b_1-a_1^Tx; \nonumber
\end{eqnarray}
or through a tight SDP relaxation \cite{J13}:
\begin{eqnarray}
({\rm P})~~~  &\min & \frac{1}{2}{\rm trace}(QX)+c^Tx \label{tight SDP}\\
&{\rm s.t.}&{\rm trace}(X)\le 1, \nonumber \\
&&a_i^Tx\le b_i,~i=1,\ldots,m,\nonumber \\
&&\left(\begin{array}{cc}X&x\\x^T&1\end{array}\right)\succeq 0.\nonumber
\end{eqnarray}
In other words, the polynomial solvability is built by way of
finding the hidden convexity from the non-convex problems (${\rm
T_{m}}$). The scheme is easily seen to be exorbitant as there are
examples in the same papers pointing out that neither the SOCP/SDP
reformulation nor the SDP relaxation is tight for general cases of
$m=2$ \cite{B13} and of $m=1$ \cite{J13}, respectively. According to
Burer and Anstreicher \cite{B13}, ``the computational complexity of
solving an extended trust region problem is highly dependent on the
geometry of the feasible set.''

Our basic idea to cope with the complication of the geometry is to
think the structure of the polytope directly and reduce the problem
(${\rm T_{m}}$) inductively until (${\rm T_{0}}$) is reached. To
avoid triviality, we assume, throughout the paper, that $Q$ has at
least one negative eigenvalue, i.e., $\lambda_{\min}(Q)<0$. Then,
the global minimum of (${\rm T_{m}}$) must lie on the boundary. The
boundary could be part of the unit sphere $x^T x =1$ or part of the
boundary of the polytope intersecting with the unit ball $x^T x \le
1$. In the former case when the global minimum of (${\rm T_{m}}$)
happens to be {\it solely} on the unit sphere (meaning that it does
not lie simultaneously on any boundary of the polytope), it must be
at least a local minimum of the trust region subproblem (${\rm
T_{0}}$). This case is polynomially checkable due to an early result
of Mart\'{\i}nez \cite{M94}. In the latter case if it lies on the
boundary of the polytope intersecting with the unit ball, it can be
found by solving one of the following $m$ subproblems:
\begin{eqnarray}
v({\rm T_m^j}):= &\min & f(x)= \frac{1}{2}x^TQx+c^Tx \nonumber\\
&{\rm s.t.}&x^Tx\le 1,\nonumber\\
&&a_j^Tx= b_j, \label{y:e1}\\
&&a_i^Tx\le b_i,~i=1,\ldots,j-1,j+1,\ldots,m, \nonumber
\end{eqnarray}
where the superscript $j$ varies from $1$ to $m$. By eliminating one
variable using (\ref{y:e1}), problem (${\rm T_m^j}$) can be reduced
to a type of problem (${\rm T_{m-1}}$) of $n-1$ dimensional. The
procedure can be inductively applied to (${\rm T_{m-1}}$), (${\rm
T_{m-2}}$), \ldots, and so forth until we run down to one of the
three possibilities: either an infeasible subproblem, or a convex
programming subproblem; or a subproblem of no linear inequality
constraint, i.e., (${\rm T_{0}}$).
Since $m$ is
fixed, the number of reduction can not grow exponentially and we
thus conclude the polynomial complexity of (${\rm T_{m}}$) for any
fixed positive integer $m$.

When $m$ is a variable, our induction argument eventually leads to
solve an exponential number of (${\rm T_{0}}$) so it still has to
face the curse of dimensionality. However, when $m$ is small, this
method can be very efficient. For example, when $m=1$, it requires
to only solve two subproblems of (${\rm T_{0}}$) (see Section 2
below). By the result of Mart\'{\i}nez \cite{M94}, (${\rm T_{0}}$)
has a spherical structure of global optimal solution set and
possesses at most one local non-global minimizer. To solve it
amounts to finding the root of a one-variable convex secular
function and hence avoid a generally more tedious large scale
SOCP/SDP.

Finally, we provide a new result which improves the dimension
condition [DC] in \cite{J13} to become
\begin{equation}
[{\rm NewDC}]~~ {\rm rank}\left([Q-\lambda_{\min}(Q)I_n~a_1~\ldots~a_m]\right)
\le n-1, \nonumber
\end{equation}
under which (${\rm T_m}$) admits an exact SDP relaxation. We use
some example to demonstrate that the improvement can be strict.

\section{Global Optimization of $({\rm T_{m}})$ and Complexity}

To solve (${\rm T_{m}}$), we begin with (${\rm T_{0}}$). Let $x^*$
be a global minimizer of (${\rm T_{m}}$) whereas $X_0^*$ is the set
of all the global minimizers of (${\rm T_{0}}$). Suppose $
X_0^*\bigcap \{x\mid a_i^Tx\le b_i,~i=1,\ldots,m\}\neq\emptyset. $
Then $v({\rm T_m})= v({\rm T_0})$ and any solution in the
intersection is a global optimal solution to $({\rm T_m})$.
Otherwise, $v({\rm T_m})> v({\rm T_0})$. In the former case, the
task is to find a common point $x^*$ from the intersection in
polynomial time. In the latter case, since we assume that the
smallest eigenvalue of $Q$ is negative, we need to examine every
piece of the boundary of ${\rm (T_{m})}$. In particular, if
$a_i^Tx^*<b_i,~i=1,\ldots,m$, then $x^*$ must be a local solution of
(${\rm T_{0}}$) residing on the sphere. Both cases require a
complete understanding about the structure of $X_0^*$ as well as the
local non-global minimizer of (${\rm T_{0}}$), which have been
studied in several important literature such as Mor\'{e} and
Sorensen \cite{M83}; Stern and Wolkowicz \cite{S94}; Mart\'{\i}nez
\cite{M94}; and Lucidi, Palagi, and Roma \cite{Lucidi98}. Below is a
brief review.

\subsection{The local and global minimizer of TRS}
Denote by
\[
(0>) \sigma_1=\ldots=\sigma_k<\sigma_{k+1}\le\ldots\le\sigma_n
\]
the eigenvalues of $Q$, $\Sigma={\rm
diag}(\sigma_1,\ldots,\sigma_n)$, $u_1,\ldots,u_n$ the corresponding
eigenvectors and $U=[u_1,\ldots,u_n]$. By introducing $y=Ux$,
$d=Uc$, we can express (${\rm T_{0}}$) in terms of $y$ as
\begin{eqnarray*}
v({\rm T_0})=  \min_{y^Ty\le 1}  \frac{1}{2}y^T\Sigma y+d^Ty. 
\end{eqnarray*}
Similarly, applying the same coordinate change to $a_i^Tx\le
b_i,~i=1,2,\ldots,m$ results in $a_i^TU^Ty\le b_i$. Denote by
${\tilde a}_i^T= a_i^TU^T$ and express
the polytope $\{x\mid a_i^Tx\le b_i,~i=1,\ldots,m\}$ in terms of $y$
as $\{y\mid {\tilde a}_i^Ty\le b_i,~i=1,\ldots,m\}$. 

Since Slater's condition is satisfied, all local minimizers of
(${\rm T_{0}}$) must satisfy the following KKT conditions associated
with a Lagrange multiplier $\mu\ge0$:
\begin{eqnarray}
&&(\Sigma+\mu I_n)y+d=0; \nonumber \\ 
&&\mu(y^Ty-1)=0; \nonumber \\ 
&&y^Ty\le 1. \nonumber 
\end{eqnarray}
By assumption, $\Sigma$ is not positive semidefinite, so any $y$ in
the interior of the unit-ball ($y^Ty<1$) cannot be a local
minimizer. Therefore, the necessary condition for local solutions
$y$ can be reduced to finding $\mu\ge0$ such that
\begin{eqnarray}
&&(\Sigma+\mu I_n)y+d=0; \label{nk1}\\
&&y^Ty= 1.\label{nk2}
\end{eqnarray}
When $\mu$ is not equal to any of the eigenvalues of $Q$, the matrix
$\Sigma+\mu I_n$ is invertible and hence $d\neq 0$. Then, $(y,\mu)$ with
$y_i=\frac{-d_i}{\sigma_i+\mu}$ is a solution to
(\ref{nk1})-(\ref{nk2}) if and only if $\mu$ is a root of the
following {\it secular function} \cite{S94}
\begin{equation}
\varphi(\mu)=\sum_{i=1}^n\frac{d_i^2}{(\sigma_i+\mu)^2}-1.\label{secular}
\end{equation}
From the first and second derivatives of $\varphi(\mu)$, we also
have
\begin{equation}
\phi'(\mu)=\sum_{i=1}^n\frac{-2d_i^2}{(\sigma_i+\mu)^3},\label{varphi-prime}
\end{equation}
and
$$\phi''(\mu)=\sum_{i=1}^n\frac{6d_i^2}{(\sigma_i+\mu)^4}>0,$$
which shows that the secular function is strictly convex. It turns out that
the global minimum and the local non-global minimum of (${\rm
T_{0}}$) can be distinguished by the position where their
corresponding Lagrange multiplier $\mu$ locates.

\begin{thm}\label{global-sol}(\cite{M83}) Let $(y^*, \mu^*)$ satisfy
$(\ref{nk1})-(\ref{nk2})$. Then, $y^*$ is a global minimum solution
to {\rm ($T_0$)} if and only if the Lagrange multiplier $\mu^*$
satisfies
\[
\mu^*\ge -\sigma_1(>0).
\]
\end{thm}

Based on Theorem \ref{global-sol}, in order to characterize the
global optimal solution set $X_0^*$ of (${\rm T_{0}}$), we only have
to investigate the secular function on the interval
$\mu\in[-\sigma_1,\infty)$. Notice that secular functions in
(\ref{secular}) depend on problem data $\sigma_i$ and $d_i$. Our
discussion below classifies and analyzes different types of secular
functions. The result shows that $X_0^*$ is either a singleton or a
$k$-dimensional sphere where $k$ is the multiplicity of the smallest
eigenvalue $\sigma_1$.

\begin{itemize}
\item[$\bullet$] Suppose $d_1^2+\ldots+d_k^2>0$. Then,
$\lim\limits_{\mu\rightarrow -\sigma_1^+}\varphi(\mu)=\infty;$
$\lim\limits_{\mu\rightarrow\infty}\varphi(\mu)=-1$ and
$\varphi(\mu)$ is strictly decreasing on $(-\sigma_1,\infty)$. Therefore,
the secular function $\varphi(\mu)$ has a unique solution
$\mu^*$ on $(-\sigma_1,\infty)$. In this case, $y^*$ defined by
\begin{equation}
y^*_i=-\frac{d_i}{\sigma_i+\mu^*}, ~i=1,\ldots,n.\label{y}
\end{equation}
is the unique global minimum solution of (${\rm T_{0}}$).
%
\item[$\bullet$] Suppose $d_1^2+\ldots+d_k^2=0$. There are two cases.
 \begin{itemize}
\item[(1)] $\mu^*>-\sigma_1$. It can happen
only when $d\neq 0$ and $\lim\limits_{\mu\rightarrow
-\sigma_1^+}\varphi(\mu)>0.$ Then, $y^*$ satisfies (\ref{y})
is the unique global minimizer.
    \item[(2)] $\mu^*=-\sigma_1$. By Theorem \ref{global-sol}, any $y^*$ satisfying
    \begin{eqnarray}
    &&(y_1^*)^2+\ldots+(y_k^*)^2=1-\sum_{i=k+1}^n\frac{d_i^2}{(\sigma_i-\sigma_1)^2},
    \label{y1}\\
    &&y^*_i=-\frac{d_i}{\sigma_i-\sigma_1},  ~i=k+1,\ldots,n \label{y2}
    \end{eqnarray}
    is a global minimizer. Namely, the global minimum
    solution set $X_0^*$ forms a $k$-dimensional sphere
    centered at
    $(0,\cdots,0,-\frac{d_{k+1}}{\sigma_{k+1}-\sigma_1},\cdots,-\frac{d_{n}}{\sigma_{n}-\sigma_1})$
    with the radius
    $\sqrt{1-\sum_{i=k+1}^n\frac{d_i^2}{(\sigma_i-\sigma_1)^2}}$.
    \end{itemize}
\end{itemize}

Secular functions also provide useful information on the local
non-global minimizer. In the next theorem, Mart\'{\i}nez \cite{M94}
showed that there is at most one local non-global minimizer
$\overline{y}$ in (${\rm T_{0}}$). The associated Lagrange
multiplier $\overline\mu$ is nonnegative and lies in
$(-\sigma_2,-\sigma_1).$ Moreover, Lucidi et al. \cite{Lucidi98}
showed that strict complementarity holds at the local non-global
minimizer.

\begin{thm}[\cite{M94,Lucidi98}]
Suppose $k\ge 2$ or $d_1=0$ when $k=1$, there is no local non-global
minimizer. Otherwise, there is at most one local non-global
minimizer $\overline{y}$ to {\rm ($T_0$)}, and the associated
Lagrange multiplier $\overline\mu$ satisfies
$\overline\mu\in(\max\{-\sigma_2,0\}, -\sigma_1)$ and
\begin{eqnarray}
&&\varphi(\overline\mu)=0,\label{phi1}\\
&&\varphi'(\overline\mu)\ge 0.\label{phi2}
\end{eqnarray}
Moreover, if $\overline\mu\in(\max\{-\sigma_2,0\}, -\sigma_1)$
(\ref{phi1}) and $\varphi'(\overline\mu)> 0$, then $\overline{y}$
defined as
\begin{equation}
\overline{y}_i= -\frac{d_i}{\sigma_i+\overline\mu}, ~i=1,\ldots,n \label{oy}
\end{equation}
is the unique local non-global minimizer.
\end{thm}

From the formula $\varphi'(\mu)$ in (\ref{varphi-prime}), there are
several types of convex secular functions on
$(-\sigma_2,-\sigma_1)$. It can be convex decreasing, for example,
when $d_1^2+\ldots+d_k^2=0$ and some $d_i\not=0,~i\ge k+1$ in which
case (${\rm T_{0}}$) can not have a local non-global minimizer; or
convex increasing, for example, when $d_1^2+\ldots+d_k^2>0$ and
$d_i=0,~i\ge k+1$; or have a global minimum on
$(-\sigma_2,-\sigma_1)$. In any case, the necessary conditions
(\ref{phi1})-(\ref{phi2}), once valid, must possess only a unique
solution $\overline{y}$ of the form (\ref{oy}) for
$\overline\mu\in(\max\{-\sigma_2,0\}, -\sigma_1)$ since
$\varphi(\mu)$ is strictly convex on $(-\sigma_2,-\sigma_1)$. That
is, $\overline{y}$ is only a candidate for the local non-global
minimizer of (${\rm T_{0}}$). It could otherwise represent a saddle
point rather than a local minimum.

%
%

\subsection{The intersection of $X_0^*$ and a polytope} In this subsection,
we are concerned with the following decision problem:
\begin{equation}\label{check-intersection}
X_0^*\bigcap \{y\mid {\tilde a}_i^Ty\le b_i,~i=1,\ldots,m\}\neq\emptyset.
\end{equation}
Since $X^*_0$ is a $k$-dimensional sphere as expressed in
(\ref{y1})-(\ref{y2}), we first reduce the $n$-dimensional polytope
$\{y\mid {\tilde a}_i^Ty\le b_i,~i=1,\ldots,m\}$ to the same $k$
dimension by fixing $y_i$ at $-\frac{d_i}{\sigma_i-\sigma_1}$ for
$i=k+1,\ldots,n$ and assume that
\begin{eqnarray*}
\{u\in R^k\mid \tilde
a_i^T\left(u^T,-\frac{d_{k+1}}{\sigma_{k+1}-\sigma_1},\cdots,
-\frac{d_n}{\sigma_n-\sigma_1}\right)^T\le b_i,~i=1,\ldots,m\}
\end{eqnarray*}
is non-empty. Otherwise, there would be no intersection between
$X^*_0$ and the polytope.

If $X_0^*$ is a singleton, it is easy to check because both sets are
convex. However, when $X_0^*$ is a nonconvex sphere, it is in
general difficult to determine whether (\ref{check-intersection}) is
true. Our procedure to answer yes/no for (\ref{check-intersection})
might depend exponentially on the number $m$ of linear constraints,
but only polynomially on the problem dimension $n$. Therefore, when
$m$ is a fixed constant, our method has polynomial complexity to
answer (\ref{check-intersection}).

To begin, let $L=\{u\in R^p\mid Hu\le g\}=\{u\in R^p\mid h_i^T u\le
g_i,~i=1,\ldots,m\}$; $B=\{u\in R^p\mid u^Tu\le r\}$ and $\partial
B=\{u\in R^p\mid  u^Tu=r\}$. That is, we conduct the analysis for
any $p$-dimensional space.

%

\begin{lem}\label{lem:0}
Let ${H}\in R^{m\times p}$ be column dependent and $g\in R^m$. Then,
the polytope $L$ is either infeasible or unbounded.
\end{lem}
\proof Let ${H}_1,\ldots,{H}_p$ be the columns of ${H}$. Since ${H}$
is column dependent, there is a nonzero $z\in R^p$ such that
${H}z=0$. If $u_0$ is feasible with ${H}u_0\le g$, then
${H}(u_0+\beta z)\le g$ for any scalar $\beta$. The polytope $L$ is
hence unbounded.
\endproof

\begin{lem}\label{lem:00}
Let ${H}\in R^{m\times p}$ be column independent and $g\in R^m$.
Assume that there is a $u_0$ satisfying ${H}u_0\le g$. Then, the
polytope $L$ is bounded if and only if the optimal value $f^*$ of
the following linear programming is nonnegative
\[
\begin{array}{ccc}f^*=&\min &e^T{H}u\\
&{\rm s.t.}& {H}u\le 0,\\
&& \|u\|_{\infty}\le 1,
\end{array}
\]
where $\|u\|_{\infty}:=\max_i|u_i|$. Moreover, if $f^*<0$, the
optimal solution $d$ to the linear programming is an extreme
direction of the unbounded polytope $L$.
\end{lem}
\proof Suppose $L$ is unbounded. It contains at least one extreme
ray, denoted by $\{u_0+\beta z\mid~\beta\ge 0\}$ where $z\neq 0$ and
$\|z\|_{\infty}\le 1$ such that ${H}(u_0+\beta z)\le g.$ This can
happen only when ${H}z\le 0$. Since $z\neq 0$ and $H$ is column
independent, we have ${H}z\neq 0$, i.e.,
$e^T{H}z<0$ and hence $f^*<0$. 

On the other hand, suppose $f^*<0$ and $d$ is optimal to the linear
programming. It implies that ${H}d\le0$ and $d\neq 0$. Consequently,
$\{u_0+\beta d\mid \beta\ge 0\}$ is contained in the polytope $L$,
which is therefore unbounded.
\endproof

\begin{lem}\label{lem:1}
Let $H\in R^{m\times p}$ and $g\in R^m$, where $m$ is fixed and $p$
is arbitrary. For any given $r>0$, it is polynomially checkable
whether $\{u\in R^p\mid Hu\le g, u^Tu=r\}$ is empty. Moreover, if
the set is nonempty, a feasible point can be found in polynomial
time.
\end{lem}
\proof Since both $L$ and $B$ are convex, we can either find a
${\hat u} \in L\bigcap B$ or conclude that $L\bigcap B=\emptyset$ in
polynomial time. For example, consider the convex program
\begin{equation}\label{distance-function}
\hat\delta=\min\limits_{\{(u,v)|u\in L,v\in B\}} \|u-v\|^2.
\end{equation}
If $\hat\delta>0,$ then $L\bigcap B=\emptyset.$ Otherwise, when
$\hat\delta=0,$ any optimal solution $({\hat u},{\hat v})$ to
(\ref{distance-function}) would imply that ${\hat u}={\hat v}$ is in
the intersection. If, furthermore, it happens that ${\hat u}
\in\partial B$, then $L\bigcap\partial B\neq \emptyset$. Otherwise,
we have ${\hat u}^{T}{\hat u}<r$.
%
%

Since $B$ is a {\it full} dimensional ball in $R^p$, the only
possibility that $\hat\delta=0$ but $L\bigcap\partial B=\emptyset$
is when $L$ is bounded and contained entirely in the interior of
$B$. By Lemmas \ref{lem:0} and \ref{lem:00}, the polytope $L$ is
bounded if and only if the columns of $H$ are linearly independent
and the linear programming in Lemma \ref{lem:00} has a nonnegative
optimal value. If $L$ is indeed bounded, $m\ge p$ and we enumerate
all the vertices of $L$.

Suppose $\widetilde{u}$ is a vertex point and assume, without loss
of generality, that $h_{i}^T\widetilde{u}=g_{i}$ for $i=1,...,r$ and
$h_i^T\widetilde{u}<g_i$ for $i=r+1,...,m$. Then we conclude that
rank$\{h_{1},\ldots,h_{r}\}=p$. If this is not true, there is a
vector $\eta\neq 0$ such that $h_i^T\eta=0$ for $i=1,...,r$. Then
both $\widetilde{u}+\epsilon \eta$ and $\widetilde{u}-\epsilon \eta$
are feasible solutions of $L$ for sufficiently small $\epsilon>0$,
contradicting to $\widetilde{u}$ being a vertex of $L$. Therefore,
to enumerate all extreme points of $L$, it is sufficient to pick all
$p$ linearly independent vectors out of $\{h_1,\ldots,h_m\}$ and
then check the feasibility. It follows that $L$ has at most $C(p,m)=
O(m^{\min\{p,m-p\}})$ vertices, denoted by $z_1,\ldots,z_t$. Since
$m$ is assumed to be fixed and $p\le m$, the number $t$ cannot
exceed a constant factor depending on m. It is also obvious that if
$z_i^Tz_i<r$ for $i=1,\ldots,t$, the polytope $L$ is in the strict
interior of $B$ and thus $L\bigcap\partial B=\emptyset$. If there
exists an index $j_0$ such that $z_{j_0}^Tz_{j_0}= r$, then
$z_{j_0}\in L\bigcap\partial B.$

Finally, if there is some index $j_0$ such that $z_{j_0}^Tz_{j_0}>
r$, since ${\hat u}^{T}{\hat u}<s$, the line segment $[\hat u,
z_{j_0}]$ must intersect $\partial B$ at one point. Similarly, when
$L$ is unbounded, solving the linear programming in Lemma
\ref{lem:00} yields an extreme direction $d$ at $\hat u$, along
which an intersection point at $L\bigcap \partial B$ can be easily
found. The proof is complete.
\endproof

\subsection{Iterative Reduction Procedure for Global Optimization}

Assume that $X_0^*\bigcap \{y\in R^k\mid \tilde
a_i^T\left(y^T,-\frac{d_{k+1}}{\sigma_{k+1}-\sigma_1},
\cdots,-\frac{d_n}{\sigma_n-\sigma_1}\right)^T\le
b_i,~i=1,\ldots,m\}=\emptyset$. That is, the global minimum of
$({\rm T_0})$ does not help solve $({\rm T_m})$ so that we have to
analyze directly the boundary of $\{x\in R^n | x^Tx\le1, a^T_i x\le
b_i,~i=1,2,\ldots,m\}$ and the local non-global minimizer of $({\rm
T_0}).$

The geometry of the boundary could be expressly complicate,
specified by one or several inequalities (linear or quadratic
inequalities) becoming active. However, if we consider the
boundaries ``one piece at a time'', the global minimizer $x^*$ of
$({\rm T_m})$ must belong to and thus globally minimize at least one
of the following candidate subproblems:
\begin{eqnarray}
v({\rm T_m^0}):= &\min & f(x)= \frac{1}{2}x^TQx+c^Tx \nonumber\\
&{\rm s.t.}~&x^Tx\le1,\nonumber\\
&&a_i^Tx< b_i,~i=1,2,\ldots,m;\label{T_m_0}
\end{eqnarray}
and for $j=1,2,\ldots,m$
\begin{eqnarray}
v({\rm T_m^j}):= &\min & f(x)= \frac{1}{2}x^TQx+c^Tx \nonumber\\
&{\rm s.t.}~&x^Tx\le 1,\nonumber\\
&&a_j^Tx= b_j, \label{y:e}\\
&&a_i^Tx\le b_i,~i=1,\ldots,j-1,j+1,\ldots,m. \nonumber
\end{eqnarray}
It is clear that $v({\rm T_m})=\min\{v({\rm T_m^0}),v({\rm
T_m^1}),\ldots,v({\rm T_m^m})\}$.

According to the previous analysis, there is at most one local
non-global minimizer of the trust region subproblem $({\rm T_0})$,
and the only legitimate candidate is
$\overline{x}_0=U^T\overline{y}$ where $\overline{y}$ is the unique
solution to (\ref{phi1})-(\ref{oy}) for some
$\overline\mu\in(\max\{-\sigma_2,0\}, -\sigma_1).$  Therefore, if
$f(\overline{x}_0)<\min\{v({\rm T_m^1}),\ldots,v({\rm T_m^m})\}$ and
$\overline{x}_0$ satisfies (\ref{T_m_0}), then
$$v({\rm T^0_m})\le f(\overline{x}_0)<\min\{v({\rm T_m^1}),\ldots,v({\rm
T_m^m})\}$$ and consequently $x^*$ must solve $({\rm T_m^0})$. Since
$x^*$ is in the interior of the polytope, and since the polytope has
no intersection with $X^*_0$, $x^*$ must be a local non-global
minimizer of $({\rm T_0})$. Since $\overline{x}_0$ is the unique
candidate, it follows that $x^*$ coincides with $\overline{x}_0$ and
$\overline{x}_0$ is indeed a local minimizer.

The above argument also implies that, if the unique candidate
$\overline{x}_0$ does not satisfy (\ref{T_m_0}), since there is no
alternative candidate, the optimal solution $x^*$ can not be found
from solving $({\rm T^0_m})$. In addition, when
$f(\overline{x}_0)\ge\min\{v({\rm T_m^1}),\ldots,v({\rm T_m^m})\}$,
$x^*$ should be retrieved from solving one of $({\rm
T_m^j}),~j=1,\ldots,m$ so that $({\rm T^0_m})$ need not be
considered either.

As a summary, we have
\begin{equation}\label{iter}
~~~~~~v({\rm T_m})=\left\{\begin{array}{ll}v({\rm T_0}),&{\rm if}~X_0^*\bigcap
\{x\mid a_i^Tx\le b_i,\forall i\}\neq\emptyset;\\
f(\overline{x}_0),& {\rm if}~ a_i^T\overline{x}_0< b_i,\forall i=1,\ldots,m~{\rm and}\\
& f(\overline{x}_0)<\min\{v({\rm T_m^1}),\ldots,v({\rm
T_m^m})\};\\
\min\{v({\rm T_m^1}),\ldots,v({\rm T_m^m})\},& {\rm o.w.}
\end{array}\right.
\end{equation}

It remains to show how to solve $({\rm T_m^j}),~j=1,\ldots,m$. Our
idea is to eliminate one variable using the equation (\ref{y:e}) and
maintains the same structure as minimizing a quadratic function over
the intersection of a ball centered at 0 with some polytope.

Let $P_j\in R^{n\times (n-1)}$ be a column-orthogonal matrix such
that $a_j^TP_j=0$. Let $z_0$ be a feasible solution to (\ref{y:e}).
Then $z_0-P_jP_j^Tz_0$ is also feasible to (\ref{y:e}). Using the
null-space representation, we have
\begin{equation}
\{x\in R^n\mid a_j^Tx= b_j\}=\{z_0-P_jP_j^Tz_0+P_jz\mid z\in R^{n-1}\} \label{redu}
\end{equation}
and
\begin{eqnarray*}
x^Tx&=&(z_0-P_jP_j^Tz_0+P_jz)^T(z_0-P_jP_j^Tz_0+P_jz)\\
&=&z_0^T(I-P_jP_j^T)(I-P_jP_j^T)z_0+2z_0^T(I-P_jP_j^T)P_jz+z^TP_j^TP_jz \\
&=&z_0^T(I-P_jP_j^T)z_0+z^Tz.
\end{eqnarray*}
Suppose $z_0^T(I-P_jP_j^T)z_0>1$. Then $x^Tx>1$ for all $z$ in the
null space of the column space $\{\alpha a_j|\alpha\in R\}$. It
indicates that $({\rm T_m^j})$ is infeasible since $\{x\mid~x^Tx\le
1\}\bigcap \{x\mid~a_j^Tx= b_j\}=\emptyset$. Otherwise, we can
equivalently express $({\rm T_m^j})$ as:
\begin{eqnarray}
v({\rm T_m^j}) = &\min & f(z_0-P_jP_j^Tz_0+P_jz)\nonumber\\ 
&{\rm s.t.}&z^Tz\le 1-z_0^T(I-P_jP_j^T)z_0,\nonumber\\
&&a_i^T(z_0-P_jP_j^Tz_0+P_jz)\le b_i,~i=1,\ldots,m,~i\neq j \label{Tmj}
\end{eqnarray}
which is again an extended trust region subproblem of $n-1$
variables equipped with $m-1$ linear inequality constraints. If the
subproblem $({\rm T_m^j})$ is a convex program, it can be globally
solved. Otherwise, it is reduced to an instance in $({\rm T_{m-1}})$
with at least one negative eigenvalue. Sometimes, more redundant
constraints can be also removed after the reduction. For example, we
can delete the $i$-th constraint in (\ref{Tmj}) provided
$a_i^TP_j=0$ and $a_i^T(z_0-P_jP_j^Tz_0)\le b_i$.

Iteratively applying (\ref{iter}), we will eventually terminate when
further reducing the subproblem causes (i) infeasibility; (ii) a
convex programming subproblem; or (iii) a classical trust region
subproblem (with no linear constraint left). Let $s$ be the smallest
number such that any $s+1$ inequalities are either row-dependent; or
row-independent but non-intersecting within the ball. By this
inductive way, there are at most $m\times(m-1)\times\cdots\times
(m-s+1)$ trust region subproblems to be solved. The special case
$s=2$ has been polynomially solved in \cite{B14} recently. Since $m$
is assumed to be fixed, the total number of reduction iterations is
bounded by a constant factor of $m$. We thus have proved that

\begin{thm}
For each fixed $m$, ${\rm (T_m)}$ is polynomially solvable.
\end{thm}

As examples, when $m=1$, the inductive procedure (\ref{iter})
requires to solve two trust region subproblems: one is $({\rm T_0})$
and the other one is reduced from $({\rm T_1^1})$. For $m=2$, the
three subproblems $({\rm T_0})$; $({\rm T_2^1})$ and $({\rm T_2^2})$
need be solved. The latter two can be further reduced to two trust
region subproblems each. In total, at most five trust region
subproblems are necessary for solving $({\rm T_2})$. Moreover, when
$m=1,~n\ge2$, the polytope is unbounded and there is no need to
enumerate the vertices in checking the decision problem
(\ref{check-intersection}) for a possible intersection. Same as
$m=2,~n\ge3$.

\section{Improved Dimension Condition for Exact SDP Relaxation}

In this section, we improve the very recent dimension condition by
Jeyakumar and Li \cite{J13} under which ${\rm (T_m)}$ admits an
exact SDP relaxation.

\subsection{Hidden Convexity of some special ${\rm (T_m)}$}
Without loss of generality, we may assume that
\begin{equation}
\{x\mid~a_i^Tx\le b_i,~i=1,\ldots,m\} {\rm ~has~ a~ strictly~ interior~ solution.}\label{as}
\end{equation}
Otherwise, there is a $j\in\{1,\ldots,m\}$ such that $a_j^Tx= b_j$
for all feasible $x$. According to (\ref{redu}), ${\rm (T_m)}$ can
be reduced to a similar problem with $m-1$ linear constraints.
Moreover, we can further assume Slater condition holds for ${\rm
(T_m)}$, i.e., ${\rm (T_m)}$ has a strictly interior solution. The
following proposition shows that the failure of Slater condition for
${\rm (T_m)}$ implies triviality.

\begin{prop} Under Assumption ({\ref{as}}),
${\rm (T_m)}$ has a unique feasible solution if and only if it has
no interior solution.
\end{prop}
\proof Obviously, when there is only a unique solution for ${\rm
(T_m)}$, it cannot be an interior point.  Now suppose ${\rm (T_m)}$
has two feasible solutions: $y\neq z$. Since $\alpha y+(1-\alpha)z$
is also feasible for any $\alpha\in[0,1]$, and both $\|y\|_2\le1,~
\|z\|_2\le1$, we can always obtain some feasible solution $w$ such
that $w^Tw<1$ and $a_i^Tw\le b_i,~i=1,2,\ldots,m.$

According to Assumption ({\ref{as}}), there is an $x^c$ such that
$a_i^Tx^c< b_i,~i=1,\ldots,m$. Then for sufficient small
$\epsilon>0$, we have $y(\epsilon ):=\epsilon x^c+(1-\epsilon)w$
satisfies $y(\epsilon )^Ty(\epsilon )<1$ and $a_i^Ty(\epsilon )<b_i$
for $i=1,\ldots,m$, which contradicts the no-interior assumption.
\endproof


Now we present the main result in this section.
\begin{thm}\label{No SDP gap}
Under the assumption
\begin{equation}
[{\rm NewDC}]~~ {\rm rank}\left([Q-\lambda_{\min}(Q)I_n~a_1~\ldots~a_m]\right)\le n-1,\label{ndc}
\end{equation}
we have $v({\rm T_m})=v({\rm P})$ where $({\rm P})$ is the standard
SDP relaxation of $({\rm T_m})$ as defined in (\ref{tight SDP}).
Besides, the condition [NewDC] (\ref{ndc}) is more general than [DC]
(\ref{dc}).
\end{thm}
\proof Consider the following convex quadratically constrained
quadratic program:
\begin{eqnarray}
({\rm T^c_m})~~~ &\min & \frac{1}{2}x^T(Q-\lambda_{\min}(Q)I_n)x+c^Tx+\frac{1}{2}\lambda_{\min}(Q) \nonumber\\
&{\rm s.t.}&x^Tx\le 1, \label{Tc1}\\
&&a_i^Tx\le b_i,~i=1,\ldots,m.\label{Tc2}
\end{eqnarray}
Since the feasible region of $({\rm T^c_m})$ is nonempty and
compact, $-\infty<v({\rm T^c_m})<+\infty$. According to Proposition
6.5.6 (\cite{B03}, page 380), there is no duality gap between $({\rm
T^c_m})$ and its Lagrangian dual
  \begin{eqnarray*}
  {\rm(D^c)}~~~ &\sup& \frac{1}{2}\lambda_{\min}(Q)-\tau-\lambda-\sum_{i=1}^m\mu_{i}b_i \\
   &{\rm s.t.}&\left(\begin{array}{cc}
  Q-\lambda_{\min}(Q)I_n+ 2\lambda I_n & c+\sum_{i=1}^m\mu_{i}a_i \nonumber\\
   c^T+\sum_{i=1}^m\mu_{i}a_i^T&2\tau\end{array}\right)\succeq 0.\\
   &&\lambda\ge0,~\mu\ge 0.\nonumber
\end{eqnarray*}
The conic dual of ${\rm (D^c)}$ is
\begin{eqnarray*}
({\rm P^c})~~~&\min & \frac{1}{2}{\rm trace}((Q-\lambda_{\min}(Q)I_n)X)+c^Tx+\frac{1}{2}\lambda_{\min}(Q) \nonumber\\
&{\rm s.t.}&{\rm trace}(X)\le 1, \\
&&a_i^Tx\le b_i,~i=1,\ldots,m,\\
&&\left(\begin{array}{cc}X&x\\x^T&1\end{array}\right)\succeq 0.
\end{eqnarray*}
Notice that Slater condition for ${\rm (T_m)}$ implies that $({\rm
P^c})$ has a strictly feasible solution. It is trivial to see that $
{\rm(D^c)}$ also has an interior feasible solution. According to the
conic duality theorem \cite{V96}, $v({\rm P^c})=v{\rm(D^c)}$ and
both optimal values are attained.

Let $x^*$ and $(\lambda^*,\mu^*)$ be the optimal solution of $({\rm
T^c_m})$ and  ${\rm(D^c)}$, respectively. Then $\lambda^*$ and
$\mu^*$ are also the corresponding Lagrangian multipliers of
(\ref{Tc1})-(\ref{Tc2}). Let $I^*$ be the index set of active linear
constraints at $x^*$, i.e., $I^*=\{i\mid~a_i^Tx^*=b_i\}$.

Suppose $\lambda^*>0$. By the complementarity, $x^{*T}x^*=1$ and
hence $v({\rm T_m})=v({\rm T^c_m})$. Now assume $\lambda^*=0$.
Corresponding to the Lagrangian multiplier $(0,\mu^*)$, all the
feasible solutions $x$ of $({\rm T^c_m})$ which satisfy the
following KKT conditions are optimal to $({\rm T^c_m})$:
\begin{equation}
\left\{\begin{array}{c}
(Q-\lambda_{\min}(Q)I_n)x+c+\sum_{i=1}^m\mu^*_ia_i=0,\\
a_i^Tx= b_i,~\forall i\in I^*.
\end{array}\right.\label{kkt}
\end{equation}
Since Assumption [{\rm NewDC}] (\ref{ndc}) implies that
\[
{\rm dim~ Ker} \left[Q-\lambda_{\min}(Q)I_n~a_1~\ldots~a_m
\right]^T = n - {\rm rank}\left([Q-\lambda_{\min}(Q)I_n~a_1~\ldots~a_m]\right)\ge 1,
\]
there exists $z\neq 0$ such that
$\left[Q-\lambda_{\min}(Q)I_n~a_1~\ldots~a_m \right]^Tz=0$. Then,
$x^*+\beta z,~\forall \beta\in R$ satisfies the KKT system
(\ref{kkt}) among which
\[
\|x^*+\beta^* z\|^2_2-1=(z^Tz) \beta^{*2}+2(x^{*T}z)\beta^*
+x^{*T}x^{*}-1=0
\]
where
\[
\beta^*:=\frac{-x^{*T}z+\sqrt{(x^{*T}z)^2+z^Tz(1-x^{*T}x^{*})}}{z^Tz}.
\]
It follows that $x^*+\beta^* z$ solves $({\rm T^c_m})$ and
consequently $v({\rm T_m})=v({\rm T^c_m})$. In other words, the
problem $({\rm T_m})$ can be equivalently solved by the convex
program $({\rm T^c_m})$ under Assumption [{\rm NewDC}].

To see the SDP relaxation $({\rm P})$ is indeed tight, we
introducing
$\widetilde{\lambda}=\lambda-\frac{1}{2}\lambda_{\min}(Q)$ and
reformulate ${\rm(D^c)}$ as
  \begin{eqnarray}
  {\rm(D^c)}~~~ &\sup&  -\tau-\widetilde{\lambda}-\sum_{i=1}^m\mu_{i}b_i \label{D_c reformulation} \\
   &{\rm s.t.}&\left(\begin{array}{cc}
  Q + 2\widetilde{\lambda} I_n & c+\sum_{i=1}^m\mu_{i}a_i \nonumber\\
   c^T+\sum_{i=1}^m\mu_{i}a_i^T&2\tau\end{array}\right)\succeq 0\\
   &&\widetilde{\lambda}\ge -\frac{1}{2}\lambda_{\min}(Q),~\mu\ge 0.\nonumber
\end{eqnarray}
Comparing (\ref{D_c reformulation}) with the Lagrangian dual problem
of ${\rm (T_m)}$:
\begin{eqnarray}
  {\rm(D)}~~~ &\sup& -\tau-\lambda-\sum_{i=1}^m\mu_{i}b_i  \nonumber\\
   &{\rm s.t.}&\left(\begin{array}{cc}
  Q+ 2\lambda I_n & c+\sum_{i=1}^m\mu_{i}a_i \nonumber\\
   c^T+\sum_{i=1}^m\mu_{i}a_i^T&2\tau\end{array}\right)\succeq 0\\
   &&\lambda\ge0,~\mu\ge 0,\nonumber
\end{eqnarray}
we find that $v{\rm (D)}\ge v{\rm(D^c)}$ since it is assumed that
$\lambda_{\min}(Q)<0$. Moreover, it can be easily verified that (D)
is the conic dual of the standard SDP relaxation (P). Since Slater
condition holds for both (P) and (D), it implies that $v{\rm
(P)}=v{\rm(D)}$.

As a summary, we have the following chain of inequalities
\[
v({\rm T_m})\ge v{\rm (P)}= v{\rm (D)}\ge v{\rm(D^c)}=v({\rm
T^c_m})=v({\rm T_m}),
\]
which proves $v({\rm T_m})=v{\rm(P)}$.

Finally, we can verify that [NewDC] (\ref{ndc}) actually improves
[DC] (\ref{dc}) by the following derivation:
\begin{eqnarray*}
&&{\rm dim~Ker}(Q-\lambda_{\min}(Q)I_n)\ge {\rm dim~span}\{a_1,\ldots,a_m\}+1\\
&\Longleftrightarrow&n-{\rm rank}(Q-\lambda_{\min}(Q)I_n)\ge {\rm dim~span}\{a_1,\ldots,a_m\}+1\\
&\Longleftrightarrow& {\rm rank}(Q-\lambda_{\min}(Q)I_n)+{\rm rank}([a_1,\ldots,a_m])\le n-1\\
&\Longrightarrow&{\rm rank}\left([Q-\lambda_{\min}(Q)I_n~a_1~\ldots~a_m]\right)\le n-1.
\end{eqnarray*}
The proof of the theorem is thus completed.
\endproof

For the special case $m=1$, the condition [NewDC] (\ref{ndc}) can be
easily satisfied except that when the smallest eigenvalue of $Q$
does not repeat (i.e., $k=1$, so that ${\rm
rank}\left(Q-\lambda_{\min}(Q)I_n\right)=n-1$) and $a_1$ happens to
be in ${\rm Ker}(Q-\lambda_{\min}(Q)I_n)$. Namely, we have
\begin{cor}
When $m=1$, $v({\rm T_1})>v({\rm P})$ happens only when the smallest
eigenvalue of $Q$ does not repeat and $a_1\in {\rm
Ker}(Q-\lambda_{\min}(Q)I_n)$.
\end{cor}

In general, in order to make the condition [NewDC] (\ref{ndc}) hold
for $m\ge2$, the smallest eigenvalue of $Q$ must repeat at least
once or $a_1,\ldots,a_m$ are in the range space of
$Q-\lambda_{min}(Q)I_n$.
 .

\subsection{Examples} The following examples illustrate the applicability of
the condition [NewDC]. Example \ref{ex1} shows that [NewDC] strictly
improves [DC]. Example \ref{ex2} shows that SDP is not tight when
[NewDC] fails. Finally, Example \ref{ex3} gives a special type of
$({\rm T_1})$ which has a tight SDP {\it without} any condition.

\begin{exam}\label{ex1}
Consider an instance of $({\rm T_1})$ where $n=2$, $b_1=0$ and
\[
Q=\left[\begin{array}{rr}-2&0\\0&2\end{array}\right],~ c=\left[\begin{array}{r}0\\ 0\end{array}\right],~
a^T_1=\left[\begin{array}{r}0~-1\end{array}\right].
\]
\end{exam}
Since $\lambda_{\min}(Q)=-2$,
$Q-\lambda_{\min}(Q)I_n=\left[\begin{array}{rr}0&0\\0&4\end{array}\right]$
and
\[
{\rm dim~Ker}(Q-\lambda_{\min}(Q)I_n)=1,~{\rm dim~span}\{a_1,\ldots,a_m\}=1,
\]
the dimension condition [DC] (\ref{dc}) fails. On the other hand,
the new dimension condition [NewDC] (\ref{ndc}) holds since
\[
{\rm rank}\left([Q-\lambda_{\min}(Q)I_n~a_1~\ldots~a_m]\right)
={\rm rank}\left(\left[\begin{array}{rrr}0&0&0\\0&4&-1\end{array}\right]\right)=1\le n-1=1.
\]
According to Theorem \ref{No SDP gap}, there is no relaxation gap
between $({\rm T_1})$ and $({\rm P})$ which is verified by
\[
v({\rm T_1})=-1=v({\rm P}).
\]

\begin{exam}\label{ex2}
Consider an instance of $({\rm T_1})$ where $n=2$, $b_1=0$ and
\[
Q=\left[\begin{array}{rr}-2&0\\0&2\end{array}\right],~ c=\left[\begin{array}{r}1\\ 0\end{array}\right],~
a_1=\left[\begin{array}{r}-1~0\end{array}\right].
\]
\end{exam}
One can check that $\lambda_{\min}(Q)=-2$,
$Q-\lambda_{\min}(Q)I_n=\left[\begin{array}{rr}0&0\\0&4\end{array}\right]$
and the new dimension condition [NewDC] (\ref{ndc}) fails since
\[
{\rm rank}\left([Q-\lambda_{\min}(Q)I_n~a_1~\ldots~a_m]\right)
={\rm rank}\left(\left[\begin{array}{rrr}0&0&-1\\0&4&0\end{array}\right]\right)=2\not\le n-1=1.
\]
In this example, the SDP relaxation is not tight because it has
\[
v({\rm T_1})=0>-1=v({\rm P}).
\]

\begin{exam}\label{ex3}
Let $c=0$, $m=1$ and $b_1=0$. That is, we consider
\begin{eqnarray}
&\min & f(x)=\frac{1}{2}x^TQx \nonumber\\
&{\rm s.t.}&x^Tx\le 1,\nonumber \\
&&a_1^Tx\le 0. \nonumber
\end{eqnarray}
\end{exam}
It is not difficult to see that, let $v_{\min}$ be the unit
eigenvector of $Q$ corresponding to $\lambda_{\min}(Q)$, then
$x^*=-{\rm sign}(a_1^Tv_{\min})v_{\min}$ where
\begin{eqnarray}{\rm
sign}(t)=\left\{\begin{array}{cc}1,&{\rm if~t>0} \nonumber \\-1,&{\rm
otherwise.}\nonumber \end{array}\right.
\end{eqnarray}
solves $({\rm T_1})$ and the optimal value is
$\frac{1}{2}\lambda_{\min}(Q).$ With simple verification, we can
also see that $(X^*,x^*)=(v_{\min}v_{\min}^T,-{\rm
sign}(a_1^Tv_{\min})v_{\min})$ solves the related SDP relaxation with
the same optimal value $\frac{1}{2}\lambda_{\min}(Q).$ In other
words,
\[
v({\rm T_1})=v({\rm P})=\frac{1}{2}\lambda_{\min}(Q)
\]
with no preliminary condition.

\section{Conclusion}

From the study, we observe that solving the extended trust region
subproblem (${\rm T_m}$) by finding the hidden convexity is a very
restrictive idea. Even for $m=1$, there are simple examples with a
positive relaxation gap (c.f. Example \ref{ex2}), which in turns
require a much more complicate SOCP/SDP formulation to catch the
hidden convexity. Contrarily, we think the problem more directly
from the structure of the polytope. If $m$ is fixed, since in the
last resort we can enumerate all the vertices which depends only on
$m$, the complexity for (${\rm T_m}$) should be only a function of
$m$ and thus be fixed too. This is indeed the case as we have
exhibited in this paper. Our scheme also has to base on the
understanding of results in 1990's for characterizing the global
minimum and the local non-global minimum of the classical trust
region subproblem. Although the results have been cited for quite a
number of times in literature, we still feel that they did not
receive noticeable attention. The induction technique reduces (${\rm
T_m}$) to a couple of small-sized trust region subproblems (${\rm
T_0}$). When $m$ is not too large, it can be very efficient as only
solving the root of a one-dimensional convex secular function is
needed. Our new dimension condition [NewDC] has its own interest
too. When the problem has high multiplicity of the smallest
eigenvalue or the coefficient vector of the linear inequality
constraints are in the space spanned by the eigenvectors
corresponding to non-minimal eigenvalues of the Hessian matrix of
the objective function, (${\rm T_m}$) can be solved directly by its
SDP reformulation.



\end{document}